\documentclass{elsart}
\usepackage[english]{babel}
\usepackage[latin1]{inputenc}
\usepackage{amsmath,amsfonts}
\usepackage{mathrsfs}
\newtheorem{Theorem}{Theorem}[section]

\newtheorem{Proposition}[Theorem]{Proposition}

\numberwithin{equation}{section}
\begin{document}
\begin{frontmatter}
\title{A Note on Algebraic Multigrid Methods for the Discrete Weighted Laplacian}
\author[authorlabel1]{Stefano Serra Capizzano}
\ead{s.serracapizzano@uninsubria.it}
\author[authorlabel2]{Cristina Tablino-Possio}
\ead{cristina.tablinopossio@unimib.it}
\address[authorlabel1]{Dipartimento di Fisica e Matematica,
Universit\`a dell'Insubria - sede di Como, via Valleggio, 11,
22100 Como, Italy.}
\address[authorlabel2]{Dipartimento di Matematica e
Applicazioni, Universit\`a di Milano Bicocca, Via Cozzi, 53, 20125
Milano, Italy.}
\begin{abstract}
In recent contributions, algebraic multigrid methods have been
designed and studied from the viewpoint of the spectral
complementarity. In this note we focus our efforts on specific
applications and, more precisely, on large linear systems arising from
the approximation of weighted Laplacian with various boundary
conditions. We adapt the multigrid idea to this specific setting and
we present and critically discuss a wide numerical experimentation
showing the potentiality of the considered approach.
\end{abstract}
\end{frontmatter}
\section{Introduction}
In the present note we test a specific application of a previously
proposed algebraic multigrid procedure \cite{ST-TR-2008}. In that
manuscript, we posed and partially answered the following question:
having at our disposal an optimal multigrid procedure for $A_n x=
b$, $\{A_n\}$ being a given sequence of Hermitian positive definite
matrices of increasing dimension, which are the minimal changes (if
any) to the procedure for maintaining the optimality for $B_n y= c$,
$\{B_n\}$ new sequence of matrices with $B_n=A_n+R_n$? \newline
Of course if there is no relation between $\{A_n\}$ and $\{B_n\}$ nothing can be
said. However, under the mild assumption that there exists a value
$\vartheta$ independent of $n$ such that $A_n\le \vartheta B_n$ and
$B_n\le M I_n$ with $M$ again independent of $n$, it has been
clearly shown that the smoothers can be simply adapted and the
prolongation and restriction operators can be substantially kept
unchanged. \newline 
The aim of this paper is to show the effectiveness of
this approach in a specific setting. 
More precisely, we consider linear systems $A_n(a)u=b$
arising from Finite Difference (FD) approximations of
\[
-\nabla(a(x)\nabla u(x)) =  f(x),\ x \in \Omega=(0,1)^d,\ d\ge 1,
\]
where $a(x)\ge a_0 >0$, $f(x)$ are given bounded functions and with Dirichlet boundary conditions (BCs).
Some remarks about the case of periodic or reflective BCs are also considered
(for a discussion on this topic see \cite{NCT,model-tau}).\newline
We recall that in the case $a(x)\equiv 1$, the matrix $A_n(1)$ is
structured, positive definite,  ill-conditioned, and an optimal
algebraic multigrid method is already available (see
\cite{AD,ADS,CSun2,FS1,FS2,Hu,CCS,Sun,mcirco,mcoseni,mcosine-vcycle})
according to different BCs.\newline
Hereafter, owing to the spectral equivalence between the matrix sequences $\{A_n(a)\}$ and
$\{A_n(1)\}$, the key idea is that the multigrid procedure just devised for $\{A_n(1)\}$
can be successfully applied to $\{A_n(a)\}$ too. \newline
More in general in \cite{ST-TR-2008}, we treated the case of
structured-plus-banded uniformly bounded Hermitian positive definite
linear systems, where the banded part $R_n$ which is added to the
structured coefficient matrix $A_n$ is not necessarily definite and not necessarily
structured. In our setting $A_n=A_n(1)$  is the structured part (it is
Toeplitz, circulant etc, according to BCs) and
$R_n=A_n(a-1)$ is the non-structured, non necessarily definite
contribution. \newline
%
%
However, while a theoretical analysis of the Two-Grid Method (TGM) for
structured+banded uniformly bounded Hermitian positive definite
linear systems has been given in \cite{ST-TR-2008}, in terms of the
algebraic multigrid theory by Ruge and St\"uben
\cite{RStub}, the corresponding  analysis for the multigrid method (MGM)
is not complete and deserves further attention. Here, for MGM
algorithm, we mean the simplest (and less expensive) version of the
large family of multigrid methods, i.e., the V-cycle procedure: for
a brief description of the TGM and of the V-cycle algorithms we
refer to Section \ref{sec:algmul}, while an extensive treatment can
be found in \cite{hack}, and especially in \cite{Oost}. \newline
Indeed, the numerics in this note suggest that the MGM is
optimal in the sense that (see \cite{AxN}) the cost of solving the linear
system (inverse problem) is proportional, by a pure constant not depending on $n$, to the
cost of the matrix-vector product (direct problem): in our case more details can be given and in fact:
\newline 
$\phantom{pA}$
\textbf{a.} the observed number of iterations is bounded by a constant independent of
the size of the algebraic problem; \newline
$\phantom{pA}$ \textbf{b.} the cost per iteration (in terms of
arithmetic operations) is just linear as the size of the algebraic
problem. \newline
Furthermore, given the spectral equivalence between $\{A_n(a)\}$, $a(x)\ge a_0 >0$, and
$\{A_n(1)\}$, a simpler numerical strategy could be used:
use $A_n(1)$ as preconditioner for $A_n(a)$ in a PCG method and solve
the linear systems with coefficient matrix $A_n(1)$ by MGM. Of course,
this approach is simpler to implement, but since several linear
systems have to be solved by MGM, the flop count can be more favorable in applying the MGM directly
instead of using it as solver for the preconditioner.
\newline
The paper is organized as follows.
In Section \ref{sec:algmul} we report the standard  TGM and MGM algorithms, together
with the reference theoretical results on the TGM optimal rate of
convergence, under some general and weak
assumptions. In Section \ref{sec:numexp} the proposed approach is applied to
the discrete weighted Laplacian and several numerical experiments are considered,
by varying the diffusion function $a(x)$ with respect to its analytical features. Finally,
Section \ref{sec:finrem} deals with further considerations
concerning future work and perspectives.
\section{Two-grid and Multigrid Method}\label{sec:algmul}
We carefully report the TGM and MGM algorithms and we describe the
theoretical ground on which we base our proposal. We start with the
simpler TGM  and then we describe the MGM and its interpretation as
stationary or multi-iterative method, see \cite{Smulti}.
\subsection{Algorithm definition}%
Let $n_0$ be a positive $d$-index, $d\ge 1$, and let $N(\cdot)$ be
an increasing function with respect to $n_0$. In devising a TGM
and a MGM for the linear system
$A_{n_0} x_{n_0}=b_{n_0}$,
where $A_{n_0} \in \mathbb{C}^{N({n_0})\times N({n_0})}$  and
$x_{n_0},b_{n_0} \in \mathbb{C}^{N({n_0})}$, the ingredients below
must be considered.\newline
Let $n_1 < n_0$ (componentwise) and let $p_{n_0}^{n_1}\in
\mathbb{C}^{N({n_0})\times N({n_1})}$ be a given full-rank matrix.
In order to simplify the notation, in the following we will refer to
any multi-index $n_s$ by means of its subscript $s$, so that, e.g.
$A_s:=A_{n_s}$, $b_s:=b_{n_s}$, $p_s^{s+1}:=p_{n_s}^{n_{s+1}}$, etc. \newline
With these notations, a class of stationary iterative methods of the
form
$x_s^{(j+1)}=V_s x_s^{(j)}+ \tilde b_s$
is also considered in such a way that ${\mathcal
Smooth}(x_s^{(j)},b_s,V_s,\nu_s)$ denotes the application of this
rule $\nu_s$ times, with $\nu_s$ positive integer number, at the
dimension corresponding to the index $s$. \newline
Thus, the solution of the
linear system $A_{n_0} x_{n_0}=b_{n_0}$ is obtained by applying
repeatedly the TGM iteration, where the $j^{\mathrm th}$ iteration
\[
{x}_0^{(j+1)} = {\mathcal TGM}({x}_0^{(j)}, {b}_0,
A_0,V_{0,\mathrm{pre}},\nu_{0,\mathrm{pre}},
V_{0,\mathrm{post}},\nu_{0,\mathrm{post}})
\]
is defined by the following algorithm \cite{hack}:
\[\renewcommand{\arraystretch}{1}
\begin{tabular}{ll}
\multicolumn{2}{c}{$y_0:={\mathcal TGM}({x}_0, {b}_0, A_0,V_{0,\mathrm{pre}},\nu_{0,\mathrm{pre}},V_{0,\mathrm{post}},\nu_{0,\mathrm{post}})$} \\
\hline \\
    \fbox{\begin{tabular}{l} $\tilde x_0:={\mathcal Smooth}(x_0,b_0,V_{0,\mathrm{pre}},\nu_{0,\mathrm{pre}})
    \hskip 0.75cm \phantom{pA}
    $ \end{tabular} }  & \small{\it Pre-smoothing iterations}\\
    \ \\
 \fbox{\begin{tabular}{l}
      $r_{0}:= b_0- A_0 \tilde x_0$\\
      $r_{1}:=(p_0^1)^H r_{0}$     \\
      $\mbox{Solve\ } A_{1}y_1=r_{1}, \textrm{ with }  A_{1}:=(p_0^1)^HA_0 p_0^1$ \\
      $\tilde y_0:=\tilde x_0+p_0^1 y_1$\\
  \end{tabular}
    } & \small{\it Exact Coarse Grid Correction}\\
    \ \\
    \fbox{\begin{tabular}{l} $y_0:={\mathcal Smooth}(\tilde y_0,b_0,V_{0,\mathrm{post}},\nu_{0,\mathrm{post}})\hskip 0.8cm \phantom{pA} $
\end{tabular}} & \small{\it Post-smoothing iterations} \\
\end{tabular}
\]
The first and last steps concern the application of
$\nu_{0,\mathrm{pre}}$ steps of the \emph{pre-smoothing} (or
\emph{intermediate}) iteration and of $\nu_{0,\mathrm{post}}$ steps
of the \emph{post-smoothing} iteration, respectively. Moreover,
the intermediate steps define the so called \emph{coarse grid
correction}, that depends on the projection operator $(p_0^1)^H$. In
such a way, the TGM iteration represents a classical stationary
iterative method whose iteration matrix is given by
\begin{equation}\label{explicit:tgm-exact}
TGM_0= V_{0,\mathrm{post}}^{\nu_{0,\mathrm{post}}}
\ 
  CGC_0
\ 
V_{0,\mathrm{pre}}^{\nu_{0,\mathrm{pre}}},
\end{equation}
where $CGC_0=   I_0-p_0^1
            \left[ (p_0^1)^H A_0 p_0^1  \right]^{-1}
   (p_0^1)^H A_0$
denotes the coarse grid correction iteration matrix. \newline
The names \emph{intermediate} and \emph{smoothing iteration} used above refer
to the multi-iterative terminology \cite{Smulti}: we say that a
method is multi-iterative if it is composed by at least two distinct
iterations. The idea is that these basic components should have
complementary spectral behaviors so that the whole procedure is
quickly convergent (for details see \cite{Smulti} and Sections 7.2
and 7.3 in \cite{mcirco}).
Notice that in the setting of Hermitian positive definite and uniformly
bounded sequences, the subspace where $A_0$ is ill-conditioned
corresponds to the subspace in which $A_0$ has small eigenvalues.
\newline
Starting from the TGM, the MGM can be introduced as follows: instead of
solving directly the linear system with coefficient matrix $A_1$, the projection strategy
is recursively  applied, so obtaining a {\em multigrid method}. \newline
Let us use the Galerkin
formulation and let $n_0>n_1>\ldots
>n_l>0$, with $l$ being the maximal number of recursive calls
and with $N(n_s)$ being the corresponding matrix sizes.
\newline
The corresponding MGM generates the $j^{\mathrm th}$ iteration
\[
{x}_0^{(j+1)} = {\mathcal MGM}(0, {x}_0^{(j)}, {b}_0, A_0,
V_{0,\mathrm{pre}},{\nu}_{0,\mathrm{pre}}, V_{0,\mathrm{post}},
{\nu}_{0,\mathrm{post}})
\]
according to the following algorithm:
\[\renewcommand{\arraystretch}{1}
\begin{tabular}{lll}
    \multicolumn{3}{c}{${y}_{s}:={\mathcal MGM}(s,{x}_{s},{b}_{s}, A_{s}, V_{{s},\mathrm{pre}},{\nu}_{s,\mathrm{pre}},V_{{s},\mathrm{post}},{\nu}_{s,
        \mathrm{post}})$} \\
        \hline  \\
\texttt{if}\quad $s=l$& \texttt{then} \\
     \ \\
& \hskip -2.2cm \fbox{\begin{tabular}{l} ${\mathcal Solve}(A_{s}{y}_{s}={b}_{s})$
      \hskip 1.7cm \phantom{pA} \end{tabular}} \  \small{\it Exact solution}\\
     \ \\
     \texttt{else} &&  \\
& \hskip -2.2cm \fbox{\begin{tabular}{l} ${\tilde{x}}_{s} := {\mathcal Smooth}\left({x}_{s}, {b}_{s}, V_{{s},
            \mathrm{pre}}, \nu_{s,\mathrm{pre}} \right) \hskip 0.9cm \phantom{pA} $ \end{tabular}} \  \small{\it Pre-smoothing iterations}\\
                             \ \\
& \hskip -2.2cm \fbox{\begin{tabular}{l}${r}_{s} := {b}_{s} - A_{s}{\tilde{x}}_{s}$  \hfill \small{\it Coarse Grid Correction}\\
                                            ${r}_{s+1}   := (p_{s}^{s+1})^H {r}_{s}$   \\
                                            ${y}_{s+1}\!\!:={\mathcal MGM}(s+1,{0}_{s+1},{b}_{s+1},A_{s+1},\!V_{{s+1},\mathrm{pre}},
                                                                {\nu}_{s+1,\mathrm{pre}},V_{{s+1},\mathrm{post}},{\nu}_{s+1,\mathrm{post}})$  \\
                                            $\tilde{y}_{s} := {\tilde{x}}_{s}+  p_{s}^{s+1}{y}_{s+1}$ \\
     \end{tabular}}
     \\
          \ \\
& \hskip -2.2cm  \fbox{
          \begin{tabular}{l}  $y_s:= {\mathcal Smooth}\left( \tilde{y}_s, b_s, V_{s,\mathrm{post}}, \nu_{s,
            \mathrm{post}} \right) \hskip 0.8cm \phantom{pA}
     $ \end{tabular}} \ \small{\it Post-smoothing iterations}\\
\end{tabular}
\]
where the matrix $A_{s+1} := (p_{s}^{s+1})^H A_{s} p_{s}^{s+1}$ is
more profitably computed in the so called \emph{pre-computing
phase}. \newline 
Since the MGM is again a linear fixed-point
method, the $j^{\mathrm th}$ iteration ${x}_0^{(j+1)}$ can be expressed as $MGM_{0}{x}_0^{(j)}$
$+(I_0-MGM_{0})A_0^{-1}{b}_0$,  where the iteration matrix
$MGM_{0}$ is recursively defined  according to the following rule
(see \cite{Oost}):
\begin{equation}\label{eq:MGM}
  \begin{array}{@{}l@{\;}c@{\;}l}
    MGM_{l} & = & O,      \\
    MGM_{s} & = & V_{s,\mathrm{post}}^{\nu_{s,\mathrm{post}}} \left[
      I_{s}\!-\! p_{s}^{s+1} \!\!\left(
      I_{s+1}\!-\!MGM_{s+1}\right) \!\!
      A_{s+1}^{-1}(p_{s}^{s+1})^H\, A_{s}
      \right]  V_{ s,\mathrm{pre}}^{ \nu_{s,\mathrm{pre}}},  \\
    &  & \qquad\qquad \qquad\qquad \qquad\qquad \qquad\qquad
         \qquad  s=0, \dots, l-1,
  \end{array}
\end{equation}
and with $MGM_s$ and $MGM_{s+1}$ denoting the iteration matrices
of the multigrid procedures at two subsequent levels. \newline
At the last recursion level $l$, the linear system is solved by a
direct method and hence it can be interpreted as an iterative
method converging in a single step: this motivates the chosen
initial condition  $MGM_{l} = O$.\newline
By comparing the TGM and MGM, we observe that  the coarse grid
correction operator $CGC_s$ is replaced by an approximation, since
the matrix $A_{s+1}^{-1}$ is approximated by
$\left(I_{s+1}-\!MGM_{s+1}\right)A_{s+1}^{-1}$ as implicitly
described in (\ref{eq:MGM}) for $s=0,\ldots,l-1$. In this way step
$4.$, at the highest level $s=0$, represents an approximation of
the exact solution of step $4.$ displayed in the TGM algorithm
(for the matrix analog compare (\ref{eq:MGM}) and
(\ref{explicit:tgm-exact})).
Finally, for $l=1$  the MGM
reduces to the TGM if ${\mathcal Solve}(A_1{y}_1={b}_1)$ is
${y}_1=A_1^{-1}{b}_1$.
\subsection{Some theoretical results on TGM convergence and optimality}%
In this paper we refer to the multigrid solution of
special linear systems of the form
\begin{equation}\label{bn-system}
B_n x=b, \quad B_n\in \mathbb{C}^{N(n)\times N(n)}, \ x,b
\in\mathbb{C}^{N(n)}
\end{equation}
with $\{B_n\}$ Hermitian positive definite uniformly bounded
matrix sequence, $n$ being a positive $d$-index, $d\ge 1$ and
$N(\cdot)$ an increasing function with respect to it.  More
precisely, we assume that there exists $\{A_n\}$ Hermitian
positive definite matrix sequence such that some order relation is
linking $\{A_n\}$ and $\{B_n\}$,
for $n$ large enough and we suppose that an optimal algebraic
multigrid method is available for the solution of the systems
\begin{equation}\label{an-system}
A_n x=b, \quad A_n\in \mathbb{C}^{N(n)\times N(n)},\ x,b
\in\mathbb{C}^{N(n)}.
\end{equation}
The underlying idea is to apply for the systems (\ref{bn-system}) the some
algebraic TGM and MGM considered for the systems
(\ref{an-system}), i.e., when considering 
the very same projectors. In fact,
the quoted choice will give rise to a relevant simplification, since
it is well-known that a very crucial role in MGM is played by the choice of projector operator.
\newline
In the algebraic multigrid theory some relevant convergence
results are due to Ruge and St\"uben \cite{RStub},
 to which we referred in order to
prove our convergence results. \newline
Hereafter, by $\|\cdot\|_2$ we denote the Euclidean norm on
$\mathbb{C}^m$ and the associated induced matrix norm over
$\mathbb{C}^{m\times m}$. If $X$ is Hermitian positive definite,
then its square root obtained via the Schur decomposition is well
defined and positive definite. As a consequence we can set
 $\|\cdot\|_{X}=\|X^{1/2}\cdot\|_{2}$ the Euclidean norm
weighted by $X$ on $\mathbb{C}^m$, and the associated induced
matrix norm. In addition, the notation $X\le Y$, with $X$ and $Y$  Hermitian
matrices, means that $Y-X$ is nonnegative definite. In addition
the sequence $\{X_n\}$, with $X_n$ Hermitian positive definite matrices, is a
uniformly bounded matrix sequence if there exists $M>0$
independent of $n$ such that $\|X_n\|_2\le M$, for $n$ large
enough.
\begin{Theorem}\emph{\cite{RStub}} \label{teo:TGMconv}
Let $A_0$ be a Hermitian positive definite matrix of size
$N(n_0)$, let $p_0^{1}\in \mathbb{C}^{N(n_0)\times N(n_1)}$, $n_0
>n_1$, be a given full-rank matrix and let $V_{0,\mathrm{post}}$ be
the post-smoothing iteration matrix. Suppose that there
exists $\alpha_\mathrm{post}>0$, independent of $n_0$, such that
for all $x\in \mathbb{C}^{N(n_0)}$
\begin{equation}\label{hp:1post-0}
\|V_{0,\mathrm{post}}x\|_{A_0}^2 \le  \|x\|_{A_0}^2-\alpha_\mathrm{post}\ \|x\|_{A_0D_0^{-1}A_0}^2,
\end{equation}
where $D_0$ is the diagonal matrix formed by the diagonal entries
of $A_0$. \newline Assume, also, that there exists $\beta>0$,
independent of $n_0$, such that for all $x\in \mathbb{C}^{N(n_0)}$
\begin{equation}\label{hp:2-0}
\min_{y\in \mathbb{C}^{N(n_{1})} } \| x -p_0^1 y \|_{D_0}^2 \le
\beta \ \| x \|_{A_0}^2.
\end{equation}
Then, $\beta \ge \alpha_\mathrm{post}$ and $\ \|TGM_0\|_{A_0} \le
    \sqrt{1-\alpha_\mathrm{post}/\beta}<1$.
\end{Theorem}
Notice that all the constants 
$\alpha_\mathrm{post}$ and $\beta$ are
required to be independent of the actual dimension in order to
ensure a TGM convergence rate independent of the size of the
algebraic problem. \newline
It is worth stressing that Theorem \ref{teo:TGMconv} still holds if the diagonal
matrix $D_0$ is replaced by any Hermitian positive matrix $X_0$
(see e.g. {\rm \cite{ADS}}). Thus, $X_0=I$ could be a
proper choice for its simplicity. \newline
Thus, by referring to the problem in \ref{bn-system} we can claim the following results.
\begin{Proposition}\emph{\cite{ST-TR-2008}}\label{prop:TGMconv-A}
Let $\{A_n\}$ be a matrix sequence with $A_n$ Hermitian positive
definite matrices and  let $p_0^1\in \mathbb{C}^{N(n_0)\times
N(n_1)}$ be a given full-rank matrix for any $n_0>0$ such that
there exists $\beta_A>0$ independent of $n_0$ so that for all
$x\in \mathbb{C}^{N(n_0)}$
\begin{equation}\label{hp:2A}
\min_{y\in\mathbb{C}^{N(n_1)} } \| x -p_0^1 y \|_{2}^2 \le \beta_A
\| x \|_{A_0}^2.
\end{equation}
Let $\{B_n\}$ be another matrix sequence, with $B_n$ Hermitian
positive definite matrices, such that $A_n\le \vartheta B_n$, for
$n$ large enough, with $\vartheta>0$ absolute constant. Then, for all $x\in
\mathbb{C}^{N(n_0)}$ and $n_0$ large enough, it also holds
$\beta_B=\beta_A \vartheta$ and
\begin{equation}\label{hp:2B}
\min_{y\in \mathbb{C}^{N(n_1)} } \| x -p_0^1 y \|_{2}^2 \le
\beta_B \| x \|_{B_0}^2.
\end{equation}
\end{Proposition}
Therefore, the convergence result in Theorem \ref{teo:TGMconv} holds
true also for the matrix sequence $\{B_n\}$, if the validity of condition (\ref{hp:1post-0})
it is also guaranteed. It is
worth stressing that  in the case of Richardson smoothers such
topic is not related to any partial ordering relation connecting
the Hermitian matrix sequences $\{A_n\}$ and $\{B_n\}$, i.e.
inequalities  (\ref{hp:1post-0}), and the corresponding for the pre-smoother case,
with $\{B_n\}$ instead of $\{A_n\}$, have to be proved
independently. %
\begin{Proposition}\emph{\cite{ST-TR-2008}}\label{prop:TGMconv-S}
Let $\{B_n\}$ be an uniformly bounded matrix sequence, with
$B_n$ Hermitian positive definite matrices.
For any $n_0>0$, let $V_{n,
\mathrm{pre}}=I_n-\omega_{\mathrm{pre}} B_n$, $V_{n,
\mathrm{post}}=I_n-\omega_{\mathrm{post}} B_n$ be the
pre-smoothing and post-smoothing iteration matrices, respectively
considered in the $TGM$ algorithm. Then, there exist
$\alpha_{B,\mathrm{pre}}$, $\alpha_{B,\mathrm{post}}>0$
independent of $n_0$ such that for all $x\in \mathbb{C}^{N(n_0)}$
\begin{equation}
\|V_{0,\mathrm{pre}} x\|_{B_0}^2 \le
\|x\|_{B_0}^2-\alpha_{B,\mathrm{pre}}
\|V_{0,\mathrm{pre}}x\|_{B_0^2}^2,  \label{hp:1preB} \\
\end{equation}
\begin{equation}
\|V_{0,\mathrm{post}} x\|_{B_0}^2  \le
\|x\|_{B_0}^2-\alpha_{B,\mathrm{post}} \|x\|_{B_0^2}^2.
\label{hp:1postB}
\end{equation}
\end{Proposition}
See Proposition 3 in \cite{AD} for the analogous claim in the case
of $\nu_{\mathrm{pre}}, \nu_{\mathrm{post}} >0$.
\par
In this way, according to the Ruge and St\"uben algebraic theory,
we have proved the TGM optimality, that is its convergence
rate independent of the size $N(n)$ of the involved algebraic
problem.
\begin{Theorem}\emph{\cite{ST-TR-2008}}\label{tgm-only-post}
Let $\{B_n\}$ be an uniformly bounded matrix sequence, with
$B_n$ Hermitian positive definite matrices. Under the same
assumptions of Propositions \emph{\ref{prop:TGMconv-A}} and
\emph{\ref{prop:TGMconv-S}} the TGM with only one step of
post-smoothing converges to the solution of $B_nx=b$ and its
convergence rate is independent of $N(n)$.
\end{Theorem}
Clearly, as just discussed in \cite{ST-TR-2008}, the TGM iteration  with both
pre-smoothing and post-smoothing is never worse than the TGM
iteration with only post-smoothing. Therefore Theorem
\ref{tgm-only-post} implies that the TGM  with both post-smoothing
and pre-smoothing has a convergence rate independent of the dimension
for systems with matrices $B_n$ under
the same assumptions as in Theorem \ref{tgm-only-post}. \newline
Furthermore, the same issues as before,
but in connection with the MGM, deserve to be discussed. First of all, we expect that a
more severe assumption between $\{A_n\}$ and $\{B_n\}$ has to
be fulfilled in order to infer the MGM optimality for $\{B_n\}$
starting from the MGM optimality for $\{A_n\}$. The reason is
that the TGM is just a special instance of the MGM when setting
$l=1$.\newline In the TGM setting we have assumed a one side
ordering relation: here the most natural step is to consider a two
side ordering relation, that is to assume that there exist
positive constants $\vartheta_1,\vartheta_2$ independent of $n$
such that $\vartheta_1 B_n\le A_n\le \vartheta_2 B_n$, for every
$n$ large enough. The above relationships simply represent the
spectral equivalence condition for sequences of Hermitian positive
definite matrices, which is plainly fulfilled in our setting whenever the
weight function is positive, well separated from zero, and bounded. \newline
In the context of the preconditioned conjugate
gradient method (see \cite{Axelsson}), it is well known that if
$\{P_n\}$ is a given sequence of optimal (i.e., spectrally
equivalent) preconditioners for $\{A_n\}$, then $\{P_n\}$ is
also a sequence of optimal preconditioners for $\{B_n\}$ (see
e.g. \cite{TCS}). The latter fact just follows from the
observation that the spectral equivalence is an {\em equivalence}
relation and hence is transitive. \newline In summary, we have
enough heuristic motivations in order to conjecture that the
spectral equivalence is the correct, sufficient assumption and, in
reality, the numerical experiments reported in Section
\ref{sec:numexp} give a support to the latter statement.
Refer to \cite{ST-TR-2008} for some further remark about this topic.
\section{Numerical Examples}\label{sec:numexp}
Hereafter, the aim relies in testing our TGM and MGM (standard
V-cycle according to Section \ref{sec:algmul}) applied to standard
FD approximations to
\begin{equation}\label{eq:modello}
-\nabla(a(x)\nabla u(x)) =  f(x),\ x \in \Omega=(0,1)^d, \ d\ge 1,
\end{equation}
with assigned BCs and for several choices examples of the diffusion coefficient $a(x)\ge a_0>0$.
\newline
The projectors are properly chosen according to the nature of
structured part, that depends on the imposed BCs.
For instance, in the case of Dirichlet BCs we split the arising FD matrix $A_n(a)$ as
\[
A_n(a)=a_{\min}\tau_n(A_n(1))+R_n(a), \
R_n(a)=A_n(a)-a_{\min}\tau_n(A_n(1)),
\]
where $\tau_n(A_n(1))$ denotes the FD matrix belonging to the $\tau$
(or DST-I) algebra \cite{BC} obtained in the case of $a(x)\equiv 1$
and $a_{\min}$ equals the minimum of $a(x)$ on $\bar{\Omega}$ in
order to guarantee the positivity of $R_n(a)$.
\newline On the other
hand, we will use, in general as first choice, the Richardson
smoothing/in\-ter\-me\-dia\-te iteration step twice in each
iteration, before and after the coarse grid correction, with
different values of the parameter $\omega$. In some cases better
results are obtained by considering the Gauss-Seidel method for the
pre-smoothing iteration. \newline
According to the algorithm in Section \ref{sec:algmul}, when
considering the TGM, the exact solution of the system is obtained by
using a direct solver in the immediately subsequent coarse grid
dimension, while, when considering the MGM, the exact solution of
the system is computed by the same direct solver, when the coarse
grid dimension equals $15^d$ (where $d=1$ for the one-level case and
$d=2$ for the two-level case). \newline In all tables we report the
numbers of iterations required for the TGM or MGM convergence,
assumed to be reached when the Euclidean norm of the relative
residual becomes less than $10^{-7}$. We point out that the CPU
times are consistent with the iteration counts. \newline Finally, we
stress that at every level (except for the coarsest) the structured
matrix parts are never formed since we need only to store the
nonzero Fourier coefficients of the generating function at every
level for matrix-vector multiplications. Thus, besides the $O(N(n))$
operations complexity of the proposed MGM both with respect to the
structured part and clearly with respect to the non-structured one,
the memory requirements of the structured part are also very low
since there are only $O(1)$ nonzero Fourier coefficients of the
generating function at every level. On the other hand, the
projections of the initial matrix correction $R_n(a)$ are stored at
each level according to standard sparse matrix techniques during the
pre-computing phase.
%
%
\subsection{Dirichlet BCs}
We begin by considering the FD approximation of (\ref{eq:modello}) with Dirichlet BCs
in the one-level setting. As already outlined, in this case the arising matrix sequence $\{A_n(a)\}$ can be
split as
\[
A_n(a)=a_{\min}\tau_n(A_n(1))+R_n(a), \
R_n(a)=A_n(a)-a_{\min}\tau_n(A_n(1)),
\]
where $\tau_n(A_n(1))$ and $a_{\min}$ are defined as before. More
precisely, $\{\tau_n(A_n(1))\}$ is the $\tau$/Toeplitz matrix
sequence generated by the function $f(t)=2-2\cos(t)$, $t\in
(0,2\pi]$ and $a_{\min}$ equals the minimum of $a(x)$ on $\bar{\Omega}$.
\newline
Let us consider $A_0(a)\in \mathbb{R}^{n_0\times n_0}$, with $1$-index
$n_0>0$ (according to the notation introduced in Section \ref{sec:algmul}, we refer to
any multi-index $n_s$ by means of its subscript $s$). Following \cite{FS1,Sun}, we denote by $T_0^1\in
\mathbb{R}^{n_0\times n_1}$, $n_0=2n_1+1$, the operator such that
\begin{equation}\label{def:tnk-tau}
  (T_0^1)_{i,j}=\left\{ \begin{array}{ll}
            1 & \ \ \mbox{\textrm for }\  i=2j, \ \ j=1,\ldots,n_1,  \\
            0 & \ \ \mbox{\textrm otherwise},
          \end{array} \right.
\end{equation}
and we define a projector $(p_0^1)^H$, $p_0^1\in
\mathbb{R}^{n_0\times n_1}$ as
\begin{equation} \label{def:pn-tau}
p_0^1=\frac{1}{\sqrt 2} P_0 T_0^1, \quad P_0=\ \mathrm{tridiag}_0 \ [1, 2, 1]=\tau_0(\tilde{f}), \ \tilde{f}(t)=2+2\cos(t).
\end{equation}
On the other hand,  for the smoothing/intermediate Richardson iterations, the parameters $\omega$ are chosen as
\[\omega_{\mathrm{pre}}  = {2}/{(\|f\|_\infty+\|R_n(a)\|_\infty)} \]
\[  \omega_{\mathrm{post}} = {1}/{(\|f\|_\infty+\|R_n(a)\|_\infty}), \]
and we set  $\nu_\mathrm{pre}=\nu_\mathrm{post}=1$.\newline
The first set of numerical tests refer to the following settings: $a(x)\equiv 1$, $a(x)=e^x$, $a(x)=e^x+1$,
(denoted in short as a1, a2, a3 respectively). \newline
In Table \ref{table:1d-DirichletBCs-TGM} we report the numbers of
iterations required for the TGM convergence, both in the case of the
Richardson pair, and of the Richardson + Gauss-Seidel pair. All
these results confirm the optimality of the proposed TGM in the
sense that the number of iterations is uniformly bounded by a
constant not depending on the size $N(n)$ indicated in the first
column. \newline
In Table \ref{table:1d-DirichletBCs-MGM} we report the some results,
but with respect to the V-cycle application. The numerics seems
allow to claim the optimality convergence property can be extended
to the MGM. \newline It is worth stressing that the difference in
considering Richardson or Gauss-Seidel in the pre-smoothing
iterations is quite negligible in the MGM case.\newline
In Table \ref{table:1d-DirichletBCs-TGM-a2k} we report a deeper analysis of
the TGM superlinear behavior in the a2 setting. More precisely, we consider the test functions
$a(x)=e^x+10^k$ with $k$ ranging from $0$ to $6$. The
convergence behavior is unaltered in the case of the Richardson +
Gauss-Seidel pair, while for increasing $k$ we observe that the
number of required iterations by considering the
Richardson+Richardson pair progressively approaches the reference a1
case. In fact as $k\rightarrow \infty$ the function $a(x)$ after a
proper scaling converges to the constant $1$.\newline
\begin{table}
\footnotesize \centering
\caption{Number of iterations required by TGM - one-level case with Dirichlet BCs} \label{table:1d-DirichletBCs-TGM}
\begin{tabular}[c]{ccc}
\renewcommand{\arraystretch}{1}
\begin{tabular}[c]{|c|c|c|c|}\hline
\multicolumn{4}{|c|}{Richardson+Richardson} \\ \hline
\multicolumn{1}{|c|}{$N(n)$} & \multicolumn{1}{|c|}{\quad a1}
                          & \multicolumn{1}{|c|}{\quad a2}
                          & \multicolumn{1}{|c|}{\quad a3}
\\
\hline
31   & \quad 2 & \quad 8 & \quad 5   \\
63   & \quad 2 & \quad 6 & \quad 4   \\
127  & \quad 2 & \quad 5 & \quad 4   \\
255  & \quad 2 & \quad 4 & \quad 4   \\
511  & \quad 2 & \quad 4 & \quad 3   \\
\hline
\end{tabular}
\quad
\begin{tabular}[c]{|c|c|c|c|}\hline
\multicolumn{4}{|c|}{Richardson+Gauss-Seidel} \\ \hline
\multicolumn{1}{|c|}{$N(n)$} & \multicolumn{1}{|c|}{\quad a1}
                          & \multicolumn{1}{|c|}{a2}
                          & \multicolumn{1}{|c|}{a3}
\\
\hline
31   & \quad 8 & \quad 8 & \quad 8  \\
63   & \quad 8 & \quad 8 & \quad 8  \\
127  & \quad 8 & \quad 8 & \quad 8  \\
255  & \quad 8 & \quad 8 & \quad 8  \\
511  & \quad 8 & \quad 8 & \quad 8  \\
\hline
\end{tabular}
\end{tabular}
\end{table}
\begin{table}
\footnotesize \centering
\caption{Number of iterations required by MGM - one-level case with Dirichlet BCs} \label{table:1d-DirichletBCs-MGM}
\begin{tabular}[c]{ccc}
\renewcommand{\arraystretch}{1}
\begin{tabular}[c]{|c|c|c|c|}\hline
\multicolumn{4}{|c|}{Richardson+Richardson} \\ \hline
\multicolumn{1}{|c|}{$N(n)$} & \multicolumn{1}{|c|}{\quad a1}
                          & \multicolumn{1}{|c|}{\quad a2}
                          & \multicolumn{1}{|c|}{\quad a3}
\\
\hline
15   & \quad 1 & \quad 1 & \quad 1   \\
31   & \quad 2 & \quad 8 & \quad 5   \\
63   & \quad 7 & \quad 7 & \quad 7   \\
127  & \quad 8 & \quad 8 & \quad 8   \\
255  & \quad 8 & \quad 8 & \quad 8   \\
511  & \quad 8 & \quad 8 & \quad 8   \\
\hline
\end{tabular}
\quad
\begin{tabular}[c]{|c|c|c|c|}\hline
\multicolumn{4}{|c|}{Richardson+Gauss-Seidel} \\ \hline
\multicolumn{1}{|c|}{$N(n)$} & \multicolumn{1}{|c|}{\quad a1}
                          & \multicolumn{1}{|c|}{\quad a2}
                          & \multicolumn{1}{|c|}{\quad a3}
\\
\hline
15   & \quad 1 & \quad 1 & \quad 1   \\
31   & \quad 8 & \quad 8 & \quad 8   \\
63   & \quad 9 & \quad 9 & \quad 9   \\
127  & \quad 9 & \quad 9 & \quad 9   \\
255  & \quad 9 & \quad 9 & \quad 9   \\
511  & \quad 9 & \quad 9 & \quad 9   \\
\hline
\end{tabular}

\end{tabular}
\end{table}
\begin{table}
\footnotesize \centering
\caption{Number of iterations required by TGM - one-level case with Dirichlet BCs} \label{table:1d-DirichletBCs-TGM-a2k}
\begin{tabular}[c]{ccc}
\renewcommand{\arraystretch}{1}
\begin{tabular}[c]{|c|c|c|c|c|c|c|c|}\hline
\multicolumn{8}{|c|}{Richardson+Richardson} \\ \hline
\multicolumn{2}{|c|}{} & \multicolumn{6}{|c|}{$a(x)=e^x+10^k$}
\\ \cline{3-8} \multicolumn{2}{|c|}{} & \multicolumn{6}{|c|}{$k$} \\
\hline
\multicolumn{1}{|c|}{$N(n)$} & \multicolumn{1}{|c|}{a1}
                             & \multicolumn{1}{|c|}{0}
                             & \multicolumn{1}{|c|}{1}
                             & \multicolumn{1}{|c|}{2}
                             & \multicolumn{1}{|c|}{3}
                             & \multicolumn{1}{|c|}{4}
                             & \multicolumn{1}{|c|}{5}
\\
\hline
31   &  2 &  5 & 4 & 3 & 3 & 3 & 2 \\
63   &  2 &  4 & 4 & 3 & 3 & 3 & 2   \\
127  &  2 &  4 & 4 & 3 & 3 & 3 & 2  \\
255  &  2 &  4 & 3 & 3 & 3 & 3 & 2  \\
511  &  2 &  3 & 3 & 3 & 3 & 2 & 2  \\
\hline
\end{tabular}
\quad
\begin{tabular}[c]{|c|c|c|c|c|c|c|c|}\hline
\multicolumn{8}{|c|}{Richardson+Gauss-Seidel} \\ \hline
\multicolumn{2}{|c|}{} & \multicolumn{6}{|c|}{$a(x)=e^x+10^k$}
\\ \cline{3-8} \multicolumn{2}{|c|}{} & \multicolumn{6}{|c|}{$k$} \\
\hline
\multicolumn{1}{|c|}{$N(n)$} & \multicolumn{1}{|c|}{a1}
                             & \multicolumn{1}{|c|}{0}
                             & \multicolumn{1}{|c|}{1}
                             & \multicolumn{1}{|c|}{2}
                             & \multicolumn{1}{|c|}{3}
                             & \multicolumn{1}{|c|}{4}
                             & \multicolumn{1}{|c|}{5}
\\
\hline
31   &  8 &  8 &  8 &  8 &  8 &  8 & 8 \\
63   &  8 &  8 &  8 &  8 &  8 &  8 & 8 \\
127  &  8 &  8 &  8 &  8 &  8 &  8 & 8 \\
255  &  8 &  8 &  8 &  8 &  8 &  8 & 8 \\
511  &  8 &  8 &  8 &  8 &  8 &  8 & 8  \\
\hline
\end{tabular}
\end{tabular}
\end{table}
The projector definition plainly extends to the two-level setting by using tensor arguments:
$(p_0^1)^H$ is constructed in such a way that
\begin{eqnarray}
 p_0^1&=&P_0 U_0^1\\
 P_0&=&\ \mathrm{tridiag}_{n_0^{(1)}} \ [1, 2,
1]\otimes \
\mathrm{tridiag}_{n_0^{(2)}} \ [1, 2, 1],\\
U_0^1 &=& T_{0}^{1}(n_0^{(1)})\otimes T_{0}^{1} (n_0^{(2)})
\end{eqnarray}
with $n_0^{(r)}=2 n_1^{(r)}+1$ and where $T_0^1(n_0^{(r)})\in
\mathbb{R}^{n_0^{(r)}\times n_1^{(r)}}$ is the one-level matrix
given in (\ref{def:tnk-tau}). \newline
The quoted choice represents the most trivial extension of the one-level
projector to the two-level setting and is also the
less expensive from a computational point of view: in fact,
$p_0^1=\tau_0((2+2\cos(t_1)(2+2\cos(t_2)))U_0^1$ equals
$[\tau_{n_0^{(1)}}(p(2+2\cos(t_1)))T_0^1(n_0^{(1)})] \otimes
[\tau_{n_0^{(2)}}(p(2+2\cos(t_2)))T_0^1(n_0^{(2)})]$.\newline
Tables \ref{table:2d-DirichletBCs-TGM} and \ref{table:2d-DirichletBCs-MGM} report the
number of iterations with the same notation as before and where we are considering the following
function tests: $a(x)\equiv 1$, $a(x)=e^{x_1+x_2}$, $a(x)=e^{x_1+x_2}+2$,
(denoted in short as a1, a2, a3, respectively). \newline
Though the convergence behavior in the case of the
Richardson+Richadson pair is quite slow, we can observe that the
number of MGM iterations required to achieve the convergence is
essentially the same as in the TGM. This phenomenon is probably due
to some inefficiency in considering the approximation $\|
R_n(a)\|_\infty$ in the tuning of the parameter
$\omega_\mathrm{pre}$ and $\omega_\mathrm{post}$. In fact, it is
enough to substitute, for instance, the pre-smoother with the
Gauss-Seidel method in order to preserve the optimality both in the
TGM and the MGM case.\newline
\begin{table}
\footnotesize \centering
\caption{Number of iterations required by TGM - two-level case with Dirichlet BCs} \label{table:2d-DirichletBCs-TGM}
\begin{tabular}[c]{ccc}
\renewcommand{\arraystretch}{1}
\begin{tabular}[c]{|c|c|c|c|}\hline
\multicolumn{4}{|c|}{Richardson+Richardson} \\ \hline
\multicolumn{1}{|c|}{$N(n)$} & \multicolumn{1}{|c|}{\quad a1}
                          & \multicolumn{1}{|c|}{\quad  a2}
                          & \multicolumn{1}{|c|}{\quad  a3}
\\
\hline
$31^2$   & \quad 16 & \quad 73 & \quad 38     \\
$63^2$   & \quad 16 & \quad 82 & \quad 41    \\
$127^2$  & \quad 16 & \quad 86 & \quad 43    \\
$255^2$  & \quad 16 & \quad 89 & \quad 44    \\
\hline
\end{tabular}
\quad
\begin{tabular}[c]{|c|c|c|c|}\hline
\multicolumn{4}{|c|}{Richardson+Gauss-Seidel} \\ \hline
\multicolumn{1}{|c|}{$N(n)$} & \multicolumn{1}{|c|}{\quad a1}
                          & \multicolumn{1}{|c|}{\quad a2}
                          & \multicolumn{1}{|c|}{\quad a3}
\\
\hline
$31^2$   & \quad 13 & \quad 14 & \quad 14 \\
$63^2$   & \quad 13 & \quad 15 & \quad 14 \\
$127^2$  & \quad 13 & \quad 15 & \quad 14 \\
$255^2$  & \quad 13 & \quad 15 & \quad 14 \\
\hline
\end{tabular}
\end{tabular}
\end{table}
\begin{table}
\footnotesize \centering
\caption{Number of iterations required by MGM - two-level case with Dirichlet BCs} \label{table:2d-DirichletBCs-MGM}
\begin{tabular}[c]{ccc}
\renewcommand{\arraystretch}{1}
\begin{tabular}[c]{|c|c|c|c|}\hline
\multicolumn{4}{|c|}{Richardson+Richardson} \\ \hline
\multicolumn{1}{|c|}{$N(n)$} & \multicolumn{1}{|c|}{\quad a1}
                          & \multicolumn{1}{|c|}{\quad a2}
                          & \multicolumn{1}{|c|}{\quad a3}
\\
\hline
$15^2$   & \quad 1  & \quad 1  & \quad 1  \\
$31^2$   & \quad 16 & \quad 73 & \quad 38    \\
$63^2$   & \quad 16 & \quad 83 & \quad 42   \\
$127^2$  & \quad 16 & \quad 88 & \quad 43    \\
$255^2$  & \quad 16 & \quad 90 & \quad 44    \\
\hline
\end{tabular}
\quad
\begin{tabular}[c]{|c|c|c|c|}\hline
\multicolumn{4}{|c|}{Richardson+Gauss-Seidel} \\ \hline
\multicolumn{1}{|c|}{$N(n)$} & \multicolumn{1}{|c|}{a1}
                          & \multicolumn{1}{|c|}{a2}
                          & \multicolumn{1}{|c|}{a3}
\\
\hline
$15^2$   & \quad 1  & \quad 1  & \quad 1  \\
$31^2$   & \quad 13 & \quad 14 & \quad 14 \\
$63^2$   & \quad 13 & \quad 15 & \quad 15 \\
$127^2$  & \quad 13 & \quad 15 & \quad 15 \\
$255^2$  & \quad 13 & \quad 15 & \quad 15 \\
\hline
\end{tabular}
\end{tabular}
\end{table}
Finally, in Table \ref{table:2d-DirichletBCs-MGM+CG}, we report the
number of iterations required by MGM, in the case of some other test
functions. More precisely, we are considering the $\mathcal{C}^1$
function $a(x,y)=e^{x+|y-1/2|^{3/2}}$, the $\mathcal{C}^0$ function
$a(x,y)=e^{x+|y-1/2|}$, and the piecewise constant function
$a(x,y)=1$ if $x,y<1/2$, $\delta$ otherwise, with
$\delta=10,100,1000$ (denoted in short as a4, a5, a6, a7, and a8,
respectively). Taking into account the previous remarks,  our
smoothing choice is represented by the Richardson+Gauss-Seidel pair.
Moreover, the CG choice is also investigated, both in connection to
the Richardson or the Gauss-Seidel smoother.\\
The MGM optimality is again observed, according to a proper choice
of the smoother pair. \newline
In conclusion  for keeping a proper optimal convergence, we can
claim that Gauss-Seidel is necessary and the best pair is with
conjugate gradient. The explanation of this behavior is again
possible in terms of multi-iterative procedures and spectral
complementarity: in fact while Richardson is effective essentially 
only in the high frequencies space, both Gauss-Seidel and CG are able to reduce the error also
in the middle frequencies and in addition they are robust with respect to the scaling produced
by the weight function $a$.
\begin{table}
\footnotesize \centering
\caption{Number of iterations required by MGM - two-level case with Dirichlet BCs ($\dag$ = 
more than $N(n)$ iterations required for convergence)
} \label{table:2d-DirichletBCs-MGM+CG}
\begin{tabular}[c]{ccc}
\renewcommand{\arraystretch}{1}
\begin{tabular}[c]{|c|c|c|c|c|c|}\hline
\multicolumn{6}{|c|}{Richardson+Gauss-Seidel} \\ \hline
\multicolumn{1}{|c|}{$N(n)$} & \multicolumn{1}{|c|}{a4}
                          & \multicolumn{1}{|c|}{a5}
                          & \multicolumn{1}{|c|}{a6}
                          & \multicolumn{1}{|c|}{a7}
                          & \multicolumn{1}{|c|}{a8}
\\
\hline
$15^2$   & 1   & 1  & 1  & 1   & 1  \\
$31^2$   & 14  & 14 & 13 & 13  & 13  \\
$63^2$   & 15  & 15 & 13 & 13  & 13  \\
$127^2$  & 15  & 15 & 14 & 14  & 14  \\
$255^2$  & 15  & 15 & 14 & 14  & 14  \\
\hline
\multicolumn{6}{c}{} \\
\hline
\multicolumn{6}{|c|}{Richardson+CG} \\ \hline
\multicolumn{1}{|c|}{$N(n)$} & \multicolumn{1}{|c|}{a4}
                          & \multicolumn{1}{|c|}{a5}
                          & \multicolumn{1}{|c|}{a6}
                          & \multicolumn{1}{|c|}{a7}
                          & \multicolumn{1}{|c|}{a8}
\\
\hline
$15^2$   & 1   & 1  & 1  & 1     & 1  \\
$31^2$   & 21  & 24 & 46 & 1472  & $\dag$  \\
$63^2$   & 26  & 28 & 59 & 1990  & $\dag$  \\
$127^2$  & 26  & 30 & 64 & 1783  & $\dag$  \\
$255^2$  & 27  & 31 & 60 & 1973  & $\dag$  \\
\hline
\multicolumn{6}{c}{} \\
\hline
\multicolumn{6}{|c|}{Gauss-Seidel+CG} \\ \hline
\multicolumn{1}{|c|}{$N(n)$} & \multicolumn{1}{|c|}{a4}
                          & \multicolumn{1}{|c|}{a5}
                          & \multicolumn{1}{|c|}{a6}
                          & \multicolumn{1}{|c|}{a7}
                          & \multicolumn{1}{|c|}{a8}
\\
\hline
$15^2$   & 1   & 1  & 1  & 1   & 1  \\
$31^2$   & 12  & 12 & 11 & 10  & 10  \\
$63^2$   & 12  & 12 & 11 & 10  & 10  \\
$127^2$  & 12  & 12 & 11 & 10  & 10  \\
$255^2$  & 12  & 12 & 11 & 10  & 10  \\
\hline
\end{tabular}
\end{tabular}
\end{table}

\subsection{Periodic and Reflective BCs}
Hereafter, we briefly address the case of periodic or reflective BCs. 
In particular we focus on the structured part of the splitting related 
to the FD discretization with respect to $a(x)\equiv 1$, since our multigrid
strategy is tuned just with respect to it. \newline
In the case of periodic BCs the obtained matrix
sequence is the one-level circulant matrix sequence $\{S_n(f)\}$ generated
by the function $f(t)=2-2\cos(t)$, $t\in (0,2\pi]$. Following
\cite{mcirco}, we consider the operator $T_0^1\in
\mathbb{R}^{n_0\times n_1}$, $n_0=2n_1$, such that
\[
  (T_0^1)_{i,j}=\left\{ \begin{array}{ll}
            1 & \ \ \mbox{\textrm for }\  i=2j-1, \ \ j=1,\ldots,n_1,  \\
            0 & \ \ \mbox{\textrm otherwise},
          \end{array} \right.
\]
and we define a projector $(p_0^1)^H$, $p_0^1\in
\mathbb{R}^{n_0\times n_1}$, as
$p_0^1= P_0 T_0^1$,  $P_0=\ {S}_0(p)$, $p(t)=2+2\cos(t)$.
Clearly, the arising matrices are
singular, so that we consider, for instance, the classical Strang correction
\cite{T-LAA-1995}
\[
\tilde S_{n_0}(f)= S_{n_0}(f)+ f\left(\frac{2\pi}{N({n_0})}\right)
\frac{ee^t}{N({n_0})},
\]
where $e$ is the vector of all ones. \newline
By using tensor arguments, our approach plainly extend to the two-level setting. \newline
When dealing with reflective BCs, the obtained
matrix sequence is the one-level  DCT III matrix sequence ${C_n(f)}_n$
generated by the function $f(t)=2-2\cos(t)$, $t\in (0,2\pi]$.
Following \cite{mcoseni}, we consider the operator $T_0^1\in
\mathbb{R}^{n_0\times n_1}$, $n_0=2n_1$, such that
\[
  (T_0^1)_{i,j}=\left\{ \begin{array}{ll}
            1 & \ \ \mbox{\textrm for }\  i\in \{2j-1,2j\}, \ \ j=1,\ldots,n_1,  \\
            0 & \ \ \mbox{\textrm otherwise},
          \end{array} \right.
\]
and we define a projector $(p_0^1)^H$, $p_0^1\in
\mathbb{R}^{n_0\times n_1}$, as
$p_0^1= P_0 T_0^1$, $P_0=C_0(p)$, $p(t)=2+2\cos(t)$.
Clearly, due to the singularity, we consider, for instance,
\[
\tilde C_{n_0}(f)=C_{n_0}(f)+f\left(\frac{\pi}{N({n_0})}\right)
\frac{ee^t}{N({n_0})}.
\]
Again, the two-level setting is treated by using tensor arguments. \newline
The  numerical tests performed in the case of periodic or reflective BCs have the same flavor 
as those previously reported in the case of Dirichlet BCs and hence we do not report them since
the observed numerical behavior gives the same information as in the case of Dirichlet BCs. 
\section{Concluding Remarks}\label{sec:finrem}
We have presented a wide numerical experimentation concerning a
multigrid technique for the discrete weighted Laplacian \underline{with various BCs}. 
In accordance with the theoretical study in
\cite{ST-TR-2008}, the choice of the smoothers can be done taking
into account the spectral complementarity, typical of any
multi-iterative procedure. In particular, we have noticed that when
the weight function $a$ adds further difficulties in the middle
frequencies (e.g., when $a$ is discontinuous), the use of pure smoothers like Richardson, reducing the
error only the high frequencies, is not sufficient. Conversely, both
CG and Gauss-Seidel work also reasonably well in the middle
frequencies (what is called the intermediate space in a
multi-iterative method) and in fact, in some cases, their use is
mandatory if we want to keep the optimality of the method, i.e., a
convergence within a given accuracy and within a number of
iterations not depending on the size of the considered algebraic problem.
%
\bibliographystyle{elsart-num-sort}
\bibliography{IMACS08biblio}
\end{document}